 \documentclass[12pt]{article}

 \usepackage{mathrsfs}
 \usepackage{amssymb}

 \newtheorem{lemma}{Lemma}[section]
 \newtheorem{theorem}{Theorem}[section]
 \newtheorem{corollary}{Corollary}[section]

 \def\blemma{\begin{lemma}\sl{}\def\elemma{\end{lemma}}}
 \def\btheorem{\begin{theorem}\sl{}\def\etheorem{\end{theorem}}}

 \def\beqlb{\begin{eqnarray}}\def\eeqlb{\end{eqnarray}}
 \def\beqnn{\begin{eqnarray*}}\def\eeqnn{\end{eqnarray*}}

 \def\mcr{\mathscr}
 \def\mbb{\mathbb}

 \def\qed{\quad$\Box$\medskip}

 \def\<{\langle}\def\>{\rangle}

 \def\d{{\mbox{\rm d}}}\def\e{{\mbox{\rm e}}}

 \def\bfP{\mbox{\boldmath $P$}}
 \def\bfQ{\mbox{\boldmath $Q$}}

 \def\itGamma{{\it\Gamma}}
 \def\itOmega{{\it\Omega}}

 \def\dfR{{\mbb R}}\def\dfN{{\mbb N}}

\begin{document}

\noindent{Published in: {\it Acta Applicandae Mathematicae} {\bf 74}
(2002), 93--112.}

\bigskip\bigskip

\centerline{\Large\bf Non-local Branching Superprocesses}
\smallskip
\centerline{\Large\bf and Some Related Models}

\bigskip
\centerline{Donald A. Dawson\footnote{ Supported by an NSERC Research Grant and a
Max Planck Award.}}

\smallskip
\centerline{School of Mathematics and Statistics, Carleton University,}

\centerline{1125 Colonel By Drive, Ottawa, Canada K1S 5B6}

\centerline{E-mail: ddawson@math.carleton.ca}

\bigskip
\centerline{Luis G. Gorostiza\footnote{ Supported by the CONACYT
(Mexico, Grant No.~37130-E).}}

\smallskip
\centerline{Departamento de Matem\'aticas,}

\centerline{Centro de Investigaci\'on y de Estudios Avanzados,}

\centerline{A.P. 14-740, 07000 M\'exico D. F., M\'exico}

\centerline{E-mail: gortega@servidor.unam.mx}

\bigskip
\centerline{Zenghu Li\footnote{ Supported by the NNSF (China, Grant No.~10131040).}}

\smallskip
\centerline{Department of Mathematics, Beijing Normal University,}

\centerline{Beijing 100875, P.R. China}

\centerline{E-mail: lizh@email.bnu.edu.cn}

\bigskip\bigskip

{\narrower{\narrower

\centerline{\bf Abstract}

\bigskip

A new formulation of non-local branching superprocesses is given
from which we derive as special cases the rebirth, the multitype,
the mass-structured, the multilevel and the
age-reproduction-structured superprocesses and the
superprocess-controlled immigration process. This unified treatment
simplifies considerably the proof of existence of the old classes of
superprocesses and also gives rise to some new ones.

\bigskip


\bigskip


\par}\par}

\bigskip\bigskip

\section{Introduction}

\setcounter{equation}{0}

Measure-valued branching processes or superprocesses constitute a
rich class of infinite dimensional processes currently under rapid
development. Such processes arose in applications as high density
limits of branching particle systems; see e.g. Dawson (1992,
1993), Dynkin (1993, 1994), Watanabe (1968). The development of
this subject has been stimulated from different subjects including
branching processes, interacting particle systems, stochastic
partial differential equations and non-linear partial differential
equations. The study of superprocesses has also led to a better
understanding of some results in those subjects. In the
literature, several different types of superprocess have been
introduced and studied. In particular, Dawson and Hochberg (1991),
Dawson et al (1990) and Wu (1994) studied multilevel branching
superprocesses, Gorostiza and Lopez-Mimbela (1990), Gorostiza and
Roelly (1991), Gorostiza et al (1992) and Li (1992b) studied
multitype superprocesses, Dynkin (1993, 1994) and Li (1992a, 1993)
studied non-local branching superprocesses, Gorostiza (1994)
studied mass-structured superprocesses, Hong and Li (1999) and Li
(2002) studied superprocess-controlled immigration processes, and
Bose and Kaj (2000) studied age-reproduction-structured
superprocesses. Those models arise in different circumstances of
application and are of their own theoretical interests.

In this paper, we provide a unified treatment of the above models.
We first give a new formulation of the non-local branching
superprocess as the high density limit of some specific branching
particle systems. Then we derive from this superprocess the
multitype, the mass-structured, the multilevel and the
age-reproduction-structured superprocesses and
superprocess-controlled immigration processes. Another related
model, the so-called rebirth superprocesses, is also introduced to
explain the non-local branching mechanism. This unified treatment
simplifies considerably the proof of existence of the old classes
of superprocesses and also gives rise to some new ones. We think
that this treatment may give  some useful perspectives for those
models. The unification is done by considering an enriched
underlying state space $E\times I$ instead of $E$. In this way,
the mutation in types of the offspring can be modeled by jumps in
the $I$-coordinates so that the multitype superprocess can be
derived. The superprocess-controlled immigration process is
actually a special form of the multitype superprocess. To get the
mass-structured superprocess we let $I = (0,\infty)$, which
represents the mass or size of the infinitesimal particles. For
the age-reproduction-structured superprocess, we take
$I=[0,\infty) \times \dfN$, where $\dfN$ is the set of
non-negative integers, to keep the information on ages and numbers
of offspring of the particles. In this model we have of course
that any particle has non-decreasing $[0,\infty)\times
\dfN$-coordinates starting from $(0,0)$ at its birth time. To get
a two level superprocess, we simply assume that $I= M(S)^\circ$ is
the space of nontrivial finite measures on another space $S$, and
the $M(S)^\circ$-coordinate of the underlying process is a
superprocess itself. For two-level branching systems, what has
been done so far for the second level branching is local, that is,
when a superparticle branches, superoffspring are produced as
exact copies of their parent. Since the superparticles have an
internal dynamics and evolve as branching systems themselves, it
is desirable to have the possibility that the superoffspring have
internal structures different from those of their parents, which
requires a non-local branching mechanism. Models of this type have
potential applications in genetics, population dynamics and other
complex multilevel systems; see e.g. Dawson (2000) and Jagers
(1995).

Notation and basic setting: Suppose that $E$ is a Lusin
topological space, i.e., a homeomorph of a Borel subset of some
compact metric space, with Borel $\sigma$-algebra ${\mcr B}(E)$.
Let $M(E)$ denote the space of finite Borel measures on $E$
topologized by the weak convergence topology, so it is also a
Lusin topological space. Let $N(E)$ be the subspace of $M(E)$
consisting of integer-valued measures on $E$ and let $M(E)^\circ =
M(E) \setminus\{0\}$, where $0$ denotes the null measure on $E$.
The unit mass concentrated at a point $x\in E$ is denoted by
$\delta_x$. Let
 \beqnn
B(E) &=& \mbox{$\{$ bounded ${\mcr B}(E)$-measurable functions on
$E$ $\}$},
 \qquad\qquad  \\
C(E) &=& \mbox{$\{$ $f$: $f\in B(E)$ is continuous $\}$},    \\
B_a(E) &=& \mbox{$\{$ $f$: $f\in B(E)$ and $\|f\|\le a$ $\}$},
 \eeqnn
where $a\ge 0$ and ``$\|\cdot\|$''denotes the supremum norm. The subsets of positive
members of the function spaces are denoted by the superscript ``+''; e.g., $B^+(E)$,
$C^+(E)$. For $f\in B(E)$ and $\mu\in M(E)$, we write $\mu(f)$ for $\int_E f\d\mu$.

\section{Non-local branching particle systems}

\setcounter{equation}{0}

Non-local branching particle systems have been considered by many
authors. We here adapt the model of Dynkin (1993). Let $\xi=
(\itOmega,\xi_t, {\mcr F}, {\mcr F}_t, \bfP_x)$ be a right
continuous strong Markov process with state space $E$ and
transition semigroup $(P_t)_{t\ge0}$. Let $\gamma \in B^+(E)$ and
let $F(x,\d\nu)$ be a Markov kernel from $E$ to $N(E)$ such that
 \beqlb\label{2.1}
\sup_{x\in E}\int_{N(E)}\nu(1) F(x,\d\nu) <\infty.
 \eeqlb
A {\it branching particle system} with parameters $(\xi,\gamma,F)$ is described by the
following properties:

(2.A) The particles in $E$ move randomly according to the law given by the
transition  probabilities of $\xi$.

(2.B) For a particle which is alive at time $r$ and follows the path $(\xi_t)_{t\ge
r}$, the conditional probability of survival during the time interval $[r,t]$ is
$\rho(r,t):=\exp \{-\int_r^t\gamma(\xi_s)\d s\}$.

(2.C) When a particle dies at a point $x\in E$, it gives birth to a random number of
offspring in $E$ according to the probability kernel $F(x,\d\nu)$. The offspring then
start to move from their locations. (Thus the name ``non-local branching'' is used.)

In the model, it is assumed that the migrations, the lifetimes and
the branchings of the particles are independent of each other. Let
$X_t(B)$ denote the number of particles in $B\in {\mcr B}(E)$ that
are alive at time $t\ge0$ and assume $X_0(E) < \infty$. Then
$\{X_t: t\ge0\}$ is a Markov process with state space $N(E)$. For
$\sigma\in N(E)$, let $\bfQ_\sigma$ denote the conditional law of
$\{X_t: t\ge0\}$ given $X_0 = \sigma$. For $f\in B^+(E)$, put
 \beqlb\label{2.2}
u_t(x) \equiv u_t(x,f) = -\log \bfQ_{\delta_x}\exp\{-X_t(f)\}.
 \eeqlb
The independence hypotheses imply that
 \beqlb\label{2.3}
\bfQ_\sigma\exp\{-X_t(f)\} = \exp\{-\sigma(u_t)\}.
 \eeqlb
Moreover, we have the following fundamental equation
 \beqlb\label{2.4}
\e^{-u_t(x)} = \bfP_x \{\rho(0,t)\e^{-f(\xi_t)}\} +
\bfP_x\bigg\{\int_0^t \bigg[\rho(0,s)\gamma(\xi_s)
\int_{N(E)}\e^{-\nu(u_{t-s})}F(\xi_s,\d\nu)\bigg]\d s\bigg\}.
 \eeqlb
This equation is obtained by thinking that if a particle starts moving
from point $x$ at time $0$, it follows a path of $\xi$ and does not branch before
time $t$, or it first splits  at time $s\in(0,t]$. By a standard argument one sees that
equation (\ref{2.4}) is equivalent to
 \beqlb\label{2.5}
\e^{-u_t(x)} &=& \bfP_x\e^{-f(\xi_t)}
- \bfP_x\bigg\{\int_0^t \gamma(\xi_s)\e^{-u_{t-s}(\xi_s)}\d s\bigg\}      \\
& &+ \,\bfP_x\bigg\{\int_0^t\bigg[\gamma(\xi_s)
\int_{N(E)}\e^{-\nu(u_{t-s})}F(\xi_s,\d\nu)\bigg]\d s\bigg\};
\nonumber
 \eeqlb
see e.g. Dawson (1992, 1993) and Dynkin (1993, 1994). It is sometimes more convenient to denote
 \beqlb\label{2.6}
v_t(x) \equiv v_t(x,f) = 1-\exp\{-u_t(x)\},
 \eeqlb
and rewrite (\ref{2.5}) into the form
 \beqlb\label{2.7}
v_t(x) &=& \bfP_x \left\{1-\e^{-f(\xi_t)}\right\}
- \bfP_x\bigg\{\int_0^t \gamma(\xi_s) v_{t-s}(\xi_s) \d s\bigg\}   \\
& & + \,\bfP_x \bigg\{\int_0^t\bigg[\gamma(\xi_s)\int_{N(E)}
(1-\e^{-\nu(u_{t-s})}) F(\xi_s,\d\nu)\bigg]\d s\bigg\}.  \nonumber
 \eeqlb

\section{Non-local branching superprocesses}

\setcounter{equation}{0}

In this section, we prove a limit theorem for a sequence of non-local branching
particle systems. Although the particle systems considered here are very specific,
they lead to the same class of non-local branching superprocesses constructed in
Dynkin (1993, 1994) and Li (1992a) with a slightly different formulation. We shall
give some details of the derivation to clarify the meaning of the parameters,
which is needed in understanding the connections of non-local branching
with other related models.

Let $\{X_t(k): t\ge0\}$, $k=1,2,\dots$ be a sequence of branching particle
systems with parameters $(\xi,\gamma_k,F_k)$. Then for each $k$,
 \beqlb\label{3.1}
\{X_t^{(k)} := k^{-1}X_t(k): t\ge0\}
 \eeqlb
defines a Markov process in $M_k(E):=\{k^{-1}\sigma: \sigma\in
N(E)\}$. For $\sigma \in M_k(E)$, let $\bfQ_\sigma^{(k)}$ denote
the conditional law of $\{X_t^{(k)}: t\ge0\}$ given
$X_0^{(k)}=\sigma$. By (\ref{2.3}) we have
 \beqlb\label{3.2}
\bfQ_\sigma^{(k)}\exp\left\{-X_t^{(k)}(f)\right\} = \exp
\left\{-\sigma(ku^{(k)}_t)\right\},
 \eeqlb
where $u^{(k)}_t(x)$ is determined by
 \beqlb\label{3.3}
v_t^{(k)}(x) = k[1-\exp\{-u_t^{(k)}(x)\}].
 \eeqlb
and
 \beqlb\label{3.4}
v^{(k)}_t(x) &=& \bfP_x \left\{k(1-\e^{-f(\xi_t)/k})\right\}
- \bfP_x\bigg\{\int_0^t \gamma_k(\xi_s) v^{(k)}_{t-s}(\xi_s) \d s \bigg\}   \\
& &+ \,\bfP_x \bigg\{\int_0^t\bigg[ k\gamma_k(\xi_s)\int_{N(E)}
(1-\e^{-\nu(u^{(k)}_{t-s})}) F_k(\xi_s,\d\nu)\bigg]\d s\bigg\}.
\nonumber
 \eeqlb
For $\mu\in M(E)$, let $\sigma_{k\mu}$ be a Poisson random measure
on $E$ with intensity $k\mu$,  and let $\bfQ_{(\mu)}^{(k)}$ denote
the conditional law of $\{X_t^{(k)}:t\ge 0\}$ given
$X_0^{(k)}=k^{-1}\sigma_{k\mu}$. From (\ref{3.2}) we get
 \beqlb\label{3.5}
\bfQ_{(\mu)}^{(k)}\exp\left\{-X_t^{(k)}(f)\right\} =
\exp\left\{-\mu(v^{(k)}_t)\right\}.
 \eeqlb

It is natural to treat separately the offspring that start their motion from the
death sites of their parents. Suppose that $g_k\in B^+(E\times [0,1])$ and, for each
$x\in E$,
 \beqnn
g_k(x,z) = \sum_{i=0}^\infty p_i^{(k)}(x) z^i,
\qquad z\in [0,1],
 \eeqnn
is a probability generating function with $\sup_{x\in E} (\d/\d z) g_k(x,1^-) <
\infty$. Let $\alpha_k$ and $\beta_k \in B^+(E)$ and assume $\gamma_k(x)
:= \alpha_k(x) + \beta_k(x)$ is strictly positive. Let $F^{(k)}(x,\d\nu)$ be
another probability kernel from $E$ to $N(E)$ satisfying (\ref{2.1}).
We may replace $F_k(x,\d\nu)$ by
 \beqlb\label{3.6}
\gamma_k(x)^{-1} \bigg[\alpha_k(x)\sum_{i=0}^\infty p_i^{(k)}(x)
F_0^{(i)}(x,\d\nu) + \beta_k(x) F^{(k)}(x,\d\nu)\bigg],
 \eeqlb
where $F_0^{(i)}(x,\d\nu)$ denotes the unit mass concentrated at $i\delta_x$.
Intuitively, as a particle splits at $x\in E$, the branching is of local type with
probability $\alpha_k(x)/\gamma_k(x)$ and is of non-local type with probability
$\beta_k(x)/\gamma_k(x)$. If it chooses the local branching type, the distribution of
the offspring number is $\{p_i^{(k)}(x)\}$. The non-local branching at $x\in E$ is
described by the kernel $F^{(k)}(x,\d\nu)$. Now (\ref{3.4}) turns into
 \beqlb\label{3.7}
v^{(k)}_t(x) &=& \bfP_x \left\{k(1-\e^{-f(\xi_t)/k})\right\}
- \bfP_x\bigg\{\int_0^t \gamma_k(\xi_s) v^{(k)}_{t-s}(\xi_s) \d s \bigg\} \nonumber \\
& & + \,\bfP_x \bigg\{\int_0^t k\alpha_k(\xi_s) [1-
g_k(\xi_s,\e^{-u^{(k)}_{t-s}(\xi_s)})]\d s\bigg\}   \\
& & +\,\bfP_x \bigg\{\int_0^t\bigg[ k\beta_k(\xi_s)\int_{N(E)}
[1-\e^{-\nu(u^{(k)}_{t-s})}] F^{(k)}(\xi_s,\d\nu)\bigg]\d
s\bigg\},  \nonumber
 \eeqlb
or equivalently
 \beqlb\label{3.8}
v_t^{(k)}(x) + \int_0^t\bfP_x[\phi_k(\xi_s,v_{t-s}^{(k)}(\xi_s)) +
\psi_k(\xi_s,v_{t-s}^{(k)})]\d s = \bfP_x k[1-\e^{-f(\xi_t)/k}],
 \eeqlb
where
 \beqlb\label{3.9}
\phi_k(x,z) = k\alpha_k(x) [g_k(x,1-z/k) - (1-z/k)]
 \eeqlb
and
 \beqlb\label{3.10}
\psi_k(x,f)
=
\beta_k(x) [f(x) - \zeta_k(x,f)],
 \eeqlb
where
 \beqlb\label{3.11}
\zeta_k(x,f)
=
\int_{N(E)} k(1-\exp\{\nu(\log(1-f/k))\}) F^{(k)}(x,\d\nu).
 \eeqlb

Let $M_0(E)$ denote the set of all Borel probability measures on $E$. Suppose
that $h_k\in B^+(E \times M_0(E)\times [0,1])$ and, for each $(x,\pi)\in E
\times M_0(E)$,
 \beqnn
h_k(x,\pi,z) = \sum_{i=0}^\infty q_i^{(k)}(x,\pi) z^i,
\qquad z\in [0,1],
 \eeqnn
is a probability generating function with $\sup_{x,\pi} (\d/\d z) h_k(x,\pi,1^-) <
\infty$. Suppose that $G(x,\d\pi)$ is a probability kernel from $E$ to $M_0(E)$.
We may consider a special form of the second term in (\ref{3.6}) by letting
 \beqlb\label{3.12}
F^{(k)}(x,\d\nu)
=
\int_{M_0(E)}\bigg[\sum_{i=0}^\infty q_i^{(k)}(x,\pi)
(l\pi)^{*i}(\d\nu) \bigg] G(x,\d\pi),
 \eeqlb
where $l\pi(\d\nu)$ denotes the image of $\pi$ under the map $y \mapsto \delta_y$ from
$E$ to $M(E)$ and $(l\pi)^{*i}$ denotes the $i$-fold convolution of $l\pi$. Now we have
 \beqlb\label{3.13}
\zeta_k(x,f)
=
\int_{M_0(E)} k[1 - h_k(x,\pi,1-\pi(f)/k)] G(x,\d\pi).
 \eeqlb
Intuitively, if a parent particle at $x\in E$ chooses non-local branching, it
first selects an offspring-location-distribution $\pi(x,\cdot)\in M_0(E)$ according to
the probability kernel $G(x,\d\pi)$, then gives birth to a random number of
offspring according to the distribution $\{q_i^{(k)}(x,\pi(x,\cdot))\}$, and those
offspring choose their locations in $E$ independently of each other according
to $\pi(x,\cdot)$. A similar non-local branching mechanism was considered in Li
(1992a, 1993).

In view of (\ref{3.5}) and (\ref{3.8}), it is natural to assume the sequences
$\{\phi_k\}$, $\{\beta_k\}$ and $\{\zeta_k\}$ to converge if one hopes to obtain
convergence of $\{X_t^{(k)}: t\ge 0\}$ to some process $\{X_t: t\ge 0\}$ as $k\to
\infty$.

\blemma\label{l3.1}
(i) Suppose that
 \beqlb\label{3.14}
\sum_{i=0}^\infty iq_i^{(k)}(x,\pi) \le 1
 \eeqlb
and that $\zeta_k (x,f) \to \zeta(x,f)$ uniformly on $E\times B^+_a(E)$ for each $a\ge
0$, then $\zeta(x,f)$ has representation
 \beqlb\label{3.15}
\zeta(x,f) = \lambda (x,f) + \int_{M(E)^\circ} (1-\e^{-\nu(f)}) \itGamma(x,\d\nu),
 \eeqlb
where $\lambda (x,\d y)$ is a bounded kernel on $E$, and $\nu(1) \itGamma (x,\d\nu)$ is a
bounded kernel from $E$ to $M(E)^\circ$ with
 \beqlb\label{3.16}
\lambda (x,1) + \int_{M(E)^\circ} \nu(1) \itGamma(x,\d\nu) \le 1.
 \eeqlb

(ii) A functional $\zeta(x,f)$ can be given by (\ref{3.15}) and (\ref{3.16}) if and
only if it has representation
 \beqlb\label{3.17}
\zeta(x,f)
=
\int_{M_0(E)}\bigg[d(x,\pi)\pi(f)
+ \int_0^\infty (1-\e^{-u\pi(f)}) n(x,\pi,\d u)\bigg] G(x,\d\pi),
 \eeqlb
where $d\in B^+(E\times M_0(E))$, $u n(x,\pi,\d u)$ is a bounded kernel from $E\times
M_0(E)$ to $(0,\infty)$ and $G(x,\d\pi)$ is a probability kernel from $E$ to $M_0(E)$
with
 \beqlb\label{3.18}
d(x,\pi) + \int_0^\infty u n(x,\pi,\d u) \le 1.
 \eeqlb

(iii) To each function $\zeta(\cdot,\cdot)$ given by (\ref{3.17}) and (\ref{3.18})
there corresponds a sequence of the form (\ref{3.13}) satisfying the requirement of
(i).
\elemma

{\it Proof.} (i) Note that $k(1-\e^{-f/k})$ converges to $f$ uniformly in $B^+_a(E)$.
Then
 \beqnn
\zeta_k(x,k(1-\e^{-f/k}))
=
\int_{N(E)} k(1-\exp\{\nu(f)/k\}) F^{(k)}(x,\d\nu)
 \eeqnn
converges to $\zeta(x,f)$ uniformly on $E\times B^+_a(E)$. It is
known that a metric $r$ can be introduced into $E$ so that $(E,r)$
becomes a compact metric space while the Borel $\sigma$-algebra
induced by $r$ coincides with ${\mcr B}(E)$; see e.g.
Parthasarathy (1967, p.14). Now $M(E)$ endowed with weak
convergence topology is a locally compact metrizable space. Let
$\bar M(E) = M(E) \cup \{\partial\}$ be the
one-point-compactification of $M(E)$. By (\ref{3.14}), $\{k\nu(1)
F^{(k)}(x,\d(k^{-1}\nu)): x\in E, k\ge1\}$ viewed as a family of
finite measures on $\bar M(E)$ is tight. Fix $x\in E$ and take
$\{n_i\} \subset \{n\}$ such that $k_i\nu(1) F_{k_i} (x,
\d(k_i^{-1}\nu))$ converges to some finite measure $G(x,\d\nu)$ on
$\bar M(E)$ as $i\to \infty$. It follows that
 \beqnn
\zeta(x,f)
=
\int_{M(E)^\circ}(1-\e^{-\nu(f)}) \nu(1)^{-1} G(x,\d\nu),
 \eeqnn
first for $f\in C^+(E,r)$ and then for all $f\in B^+(E)$. Now (\ref{3.15}) follows by a
simple change of the measure and (\ref{3.16}) follows from (\ref{3.14}). (ii) is
immediate. To get (iii) we may set
 \beqnn
h_k(x,\pi,z)
=
1 + d(x,\pi)(z-1)
+ k^{-1} \int_0^\infty (\e^{ku(z-1)}-1) n(x,\pi,\d u).
 \eeqnn
Observe that
 \beqnn
\frac{\d^i}{\d z^i}h_k(x,\pi,0)
\ge0,
\quad i=1,2,\dots,
 \eeqnn
and (\ref{3.14}) assures that $h_k(x,\pi,0)\ge0$. Thus for fixed $(x,a)\in E\times M_0(E)$,
$h_k(x,\pi,\cdot)$ is a probability generating function. Then we define $\zeta_k(x,f)$ by
(\ref{3.13}) so that $\zeta_k(x,f) = \zeta(x,f)$ for $(x,f)\in E\times B^+_{1/k}(E)$. \qed

\blemma\label{l3.2} {\rm (Li, 1992c)}
(i) Suppose that, for each $l\ge0$, the sequence $\phi_k(x,z)$ is uniformly Lipschitz
in $z$ on the set $E\times [0,l]$ and that $\phi_k(x,z)$ converges to some $\phi(x,z)$
uniformly as $k\to\infty$, then $\phi(x,z)$ has the representation
 \beqlb\label{3.19}
\phi(x,z) = b(x)z +
c(x)z^2 + \int_0^\infty(\e^{-zu}-1+zu) m(x,\d u),
\quad x\in E, z\ge0,
 \eeqlb
where $b\in B(E)$, $c\in B^+(E)$ and $(u\land u^2) m(x,\d u)$ is a bounded kernel from
$E$ to $(0,\infty)$.

(ii) To each function $\phi(\cdot,\cdot)$ given by (\ref{3.19}) there corresponds a
sequence of the form (\ref{3.9}) satisfying the requirement of (i).
\elemma

Based on Lemmas \ref{l3.1} and \ref{l3.2}, the following result can be proved
similarly as in Dawson (1992, 1993), Dynkin (1993, 1994) and Li (1992a, c).

\blemma\label{l3.3}
If the conditions of Lemma \ref{l3.1} (i) and Lemma \ref{l3.2} (i) are fulfilled and if
$\beta_k \to \beta\in B^+(E)$ uniformly as $k\to\infty$, then for each $a\ge0$ both
$v_t^{(k)}(x,f)$ and $ku_t^{(k)}(x,f)$ converge boundedly and uniformly on the set
$[0,a] \times E\times B^+_a(E)$ of $(t,x,f)$ to the unique bounded positive solution
$V_tf(x)$ to the evolution equation
 \beqlb\label{3.20}
V_tf(x)
+ \int_0^t\bigg\{\int_E [\phi(y,V_{t-s}f(y))
+ \psi(y,V_{t-s}f)] P_s(x,\d y)\bigg\}\d s
= P_t f(x),
\quad t\ge0,
 \eeqlb
where
 \beqlb\label{3.21}
\psi(x,f)
=
\beta(x) [f(x) - \zeta(x,f)],
\qquad x\in E, f\in B^+(E).
 \eeqlb
\elemma

By Lemma \ref{l3.3} and Dawson (1993, p.42),
 \beqlb\label{3.22}
\int_{M(E)}\e^{-\nu(f)}Q_t(\mu, \d\nu)
=
\exp\{-\mu(V_tf)\},
\qquad f\in B^+(E),
 \eeqlb
defines a transition semigroup $(Q_t)_{t\ge0}$ on $M(E)$. A Markov process $\{X_t:
t\ge0\}$ with state space $M(E)$ is called a {\it non-local branching superprocess}
with parameters $(\xi,\phi,\psi)$ if it has transition semigroup $(Q_t)_{t\ge0}$.
Condition (\ref{3.14}) means that the corresponding branching particle system has
subcritical non-local branching. In terms of the limiting superprocess, this
condition is expressed as (\ref{3.16}), which is of course a restriction of the
class of $\zeta(\cdot,\cdot)$ given by (\ref{3.15}). However, since
$\phi(x,z)+\beta(x)z$ belongs to the class defined by (\ref{3.19}), and since
$\beta\in B^+(E)$ is arbitrary, (\ref{3.16}) does not put any restriction on
the generality of
 \beqnn
\phi(x,f(x)) + \psi(x,f)
=
\phi(x,f(x)) + \beta(x)f(x) - \beta(x)\zeta(x,f).
 \eeqnn
Therefore, the class of non-local branching superprocesses given by (\ref{3.20})
and (\ref{3.22}) coincides with those constructed Dynkin (1993, 1994) and Li
(1992a), where the first term of $\psi(x,\cdot)$ was written into $\phi(x,\cdot)$.
In principle, (\ref{3.20}) and (\ref{3.22}) give the most general non-local branching
superprocesses constructed in the literature up to now. A more general class of
non-local branching superprocesses were discussed in Dynkin et al (1994), but their
existence has not been established. The next theorem follows similarly as in Li
(1992a, c).

\btheorem\label{t3.1}
Let $\{X_t^{(k)}: t\ge0\}$ be the sequence of renormalized branching particle systems
determined by (\ref{3.5}) and (\ref{3.8}), and let $\{X_t: t\ge0\}$ be the non-local
branching superprocess with transition semigroup $(Q_t)_{t\ge0}$ given by (\ref{3.20})
and (\ref{3.22}). Assume that the conditions of Lemma \ref{l3.1} (i) and Lemma
\ref{l3.2} are fulfilled. Then for every $\mu\in M(E)$, $0\le t_1<\dots <t_n$ and
$a\ge 0$, as $k\to\infty$,
 \beqnn
\bfQ_{(\mu)}^{(k)} \exp \bigg\{-\sum_{i=1}^n
X_{t_i}^{(k)}(f_i)\bigg\} \ \to \
\bfQ_{\mu}\exp\bigg\{-\sum_{i=1}^n X_{t_i}(f_i)\bigg\}
 \eeqnn
uniformly on $f_1,\dots,f_n\in B^+_a(E)$.
\etheorem

Naturally, we may regard $\{X_t^{(k)}: t\ge0\}$ as a process in
the space $M(E)$. Then the above theorem shows that the finite
dimensional distributions of $\{X_t^{(k)}: t\ge0\}$ under
$\bfQ_{(\mu)}^{(k)}$ converge as $k\to\infty$ to those of $\{X_t:
t\ge0\}$ under $\bfQ_\mu$. Therefore, the non-local branching
superprocess is a small particle approximation for the non-local
branching particle  system. Heuristically, $\xi$ gives the law of
the migration of the ``particles'', $\phi(x,\cdot)$ describes the
amount of offspring born at $x\in E$ by a parent that dies at this
point, and $\psi(x,\cdot)$ describes the amount of the offspring
born by this parent that are displaced randomly into the space
according to distributions $\pi$ randomly chosen by $G(x,d\pi)$.
Thus the locations of non-locally displaced offspring involve two
sources of randomness.

Replacing $f$ in (\ref{3.21}) and (\ref{3.22}) by $\theta f$ and differentiating at
$\theta =0$ we see that the first moments of the superprocess are given by
 \beqlb\label{3.23}
\int_{M(E)} \nu(f) Q_t(\mu, \d\nu)
=
\mu(T_tf),
\qquad t\ge0,f\in B^+(E),
 \eeqlb
where $(T_t)_{t\ge0}$ is a locally bounded semigroup of kernels on $E$
determined by
 \beqlb\label{3.24}
T_tf(x)
+ \int_0^t\bigg\{\int_E b(y)T_{t-s}f(y) + \beta(y)\big[T_{t-s}f(y)
- m(y,T_{t-s}f)\big]P_s(x,\d y)\bigg\}\d s
= P_t f(x),
 \eeqlb
and $m(x,\d y)$ is the bounded kernel on $E$ defined by
 \beqlb\label{3.25}
m(x,f) = \lambda (x,f) + \int_{M(E)^\circ} \nu(f) \itGamma(x,\d\nu).
 \eeqlb
In particular, we may define another locally bounded semigroup of kernels
$(U_t)_{t\ge0}$ on $E$ by
 \beqlb\label{3.26}
U_tf(x)
+ \int_0^t\bigg\{\int_E \beta(y)\big[U_{t-s}f(y)
- m(y,U_{t-s}f)\big]P_s(x,\d y)\bigg\}\d s
= P_t f(x),
 \eeqlb
which has weak generator $G$ such that
 \beqlb\label{3.27}
Gf(x)
=
Af(x) + \beta(x)[m(x,f) - f(x)],
\qquad f\in{\cal D}(A),
 \eeqlb
where $A$ denotes the weak generator of $(P_t)_{t\ge0}$. Now the generator $B$ of
$(T_t)_{t\ge0}$ can be expressed as
 \beqlb\label{3.28}
Bf(x)
=
Gf(x) - b(x)f(x),
\qquad f\in{\cal D}(A).
 \eeqlb
By a comparison theorem we have $T_tf\le \e^{\|b\|t} U_tf$ for all $t\ge0$ and
$f\in B(E)^+$. From this and (\ref{3.23}) we have
 \beqlb\label{3.29}
\int_{M(E)} \nu(f) Q_t(\mu, \d\nu)
\le
\e^{\|b\|t}\mu(U_tf),
\qquad t\ge0,f\in B^+(E).
 \eeqlb
Note that (\ref{3.29}) implies that $\nu\mapsto \nu(1)$ is a $\|b\|$-excessive
function for $(Q_t)_{t\ge0}$.

To conclude this section, let us consider briefly the special, and possibly more
desirable, case where $G(x,\d\pi)\equiv$ unit mass at some $\pi(x,\cdot)\in M_0(E)$,
that is, the non-locally displaced offspring born at $x\in E$ choose their locations
independently according to the (non-random) distribution $\pi(x,\cdot)$. In this case,
the non-local branching mechanism is given by
 \beqlb\label{3.30}
\psi(x,f) = \beta(x)[f(x) - \zeta (x,\pi(x,f))],
\qquad x\in E, f\in B^+(E),
 \eeqlb
where
 \beqlb\label{3.31}
\zeta (x,z) = d(x)z + \int_0^\infty (1-e^{-zu}) n(x,\d u),
\qquad x\in E,z\ge0,
 \eeqlb
where $d\in B^+(E)$ and $un(x,\d u)$ is a bounded kernel from $E$ to $(0,\infty)$ with
 \beqlb\label{3.32}
m(x) := d(x) + \int_0^\infty u n(x,\d u) \le 1,
\qquad x\in E.
 \eeqlb
In particular, if $\zeta(x,z)\equiv z$, we may rewrite (\ref{3.20}) formally as
 \beqnn
\frac{\d}{\d t}V_tf(x)
=
AV_tf(x) - \phi(x,V_tf(x)) + \beta(x)[\pi(x,V_tf)-V_tf(x)],
\quad t\ge 0, x\in E,
 \eeqnn
with initial condition $V_0f=f$. This equation corresponds to a superprocess
with underlying generator $A$ and non-trivial local and non-local branching
mechanisms. Alternatively, we may also think that the superprocess has underlying
generator $Af(x)+\beta(x)[\pi(x,f)-f(x)]$ and only non-trivial local branching
mechanism. Since the generator $B$ of a general Markov process
in $E$ is the limit of a sequence of operators of the type $\beta(x)[\pi(x,f)-f(x)]$,
in principle a superprocess with more general underlying generator $A+B$
and only local branching mechanism can be approximated by a sequence of
superprocesses with underlying generator $A$ and non-trivial local and non-local
branching mechanisms. Under suitable conditions it is also possible to establish
convergence of branching particle systems with underlying generator $A$ and non-trivial
local and non-local branching mechanisms to the superprocess with underlying
generator $A+B$ and with only non-trivial local branching mechanism, which has
been done in a particular setting in Gorostiza (1994); see also section 6.

\section{Rebirth superprocesses}

\setcounter{equation}{0}

We may consider a modification of the branching particle system described in the last
two sections. Let $(\xi,\gamma,F)$ be given as in section 2. A {\it rebirth branching
particle system} with parameters $(\xi,\gamma,F)$ is described by (2.A), (2.B) and
the following

(4.C) When a particle dies at a point $x\in E$, it gives birth to a random number of
offspring in $E$ according to the probability kernel $F(x,\d\nu)$. In addition, the
parent particle itself is replaced by an extra offspring at site $x\in E$, that is,
the parent particle is reborn. All the offspring then start to move from their
locations.

Let $\{X_t: t\ge0\}$ be the process defined in the same way as in section 2. Then
$\{X_t: t\ge0\}$ is still a Markov process with state space $N(E)$. We also have
(\ref{2.2}) and (\ref{2.3}), but (\ref{2.4}) is now replaced by
 \beqlb\label{4.1}
\e^{-u_t(x)}
&=& \bfP_x \{\rho(0,t)\e^{-f(\xi_t)}\}    \\
& & + \,\bfP_x\bigg\{\int_0^t \bigg[\rho(0,s) \gamma(\xi_s)
\int_{N(E)} \e^{-u_{t-s}(\xi_s)} \e^{-\nu(u_{t-s})}
F(\xi_s,\d\nu)\bigg]\d s\bigg\}. \nonumber
 \eeqlb
This is equivalent to
 \beqnn
\e^{-u_t(x)} &=& \bfP_x\e^{-f(\xi_t)}
- \bfP_x\bigg\{\int_0^t \gamma(\xi_s)\e^{-u_{t-s}(\xi_s)}\d s\bigg\}      \\
& &
+\,\bfP_x\bigg\{\int_0^t\bigg[\gamma(\xi_s)\int_{N(E)}\e^{-u_{t-s}(\xi_s)}
\e^{-\nu(u_{t-s})}F(\xi_s,\d\nu)\bigg]\d s\bigg\},  \nonumber
 \eeqnn
or
 \beqnn
1-\e^{-u_t(x)} &=& \bfP_x \left\{1-\e^{-f(\xi_t)}\right\}
- \bfP_x \bigg\{\int_0^t\gamma(\xi_s)(1-\e^{-u_{t-s}(\xi_s)})\d s\bigg\} \nonumber \\
& & +\,\bfP_x \bigg\{\int_0^t\bigg[\gamma(\xi_s)\int_{N(E)}
(1-\e^{-u_{t-s}(\xi_s)})\e^{-\nu(u_{t-s})} F(\xi_s,\d\nu)\bigg]\d s\bigg\}    \\
& & +\,\bfP_x \bigg\{\int_0^t\bigg[\gamma(\xi_s)\int_{N(E)}
(1-\e^{-\nu(u_{t-s})}) F(\xi_s,\d\nu)\bigg]\d s\bigg\}.  \nonumber
 \eeqnn
Let $v_t(x) \equiv v_t(x,f)$ be defined by (\ref{2.6}). Then we have
 \beqlb\label{4.2}
v_t(x) &=& \bfP_x \left\{1-\e^{-f(\xi_t)}\right\}
- \bfP_x \bigg\{\int_0^t\gamma(\xi_s) v_{t-s}(\xi_s) \d s \bigg\}  \nonumber  \\
& & +\,\bfP_x \bigg\{\int_0^t\bigg[\gamma(\xi_s)\int_{N(E)}
v_{t-s}(\xi_s) \e^{-\nu(u_{t-s})} F(\xi_s,\d\nu)\bigg]\d s\bigg\}    \\
& & +\,\bfP_x \bigg\{\int_0^t\bigg[\gamma(\xi_s)\int_{N(E)}
(1-\e^{-\nu(u_{t-s})}) F(\xi_s,\d\nu)\bigg]\d s\bigg\}.  \nonumber
 \eeqlb

We now consider a sequence of rebirth branching particle systems $\{X_t(k):
t\ge0\}$ with parameters $(\xi,\gamma_k,F_k)$. Define $\{X_t^{(k)}: t\ge0\}$ and choose
$F_k$ as in section 3 with $\alpha_k(x) \equiv 0$ and $\gamma_k(x) \equiv \beta_k(x)$.
Then (\ref{3.5}) remains valid if we replace (\ref{3.7}) by
 \beqlb\label{4.3}
v^{(k)}_t(x) &=& \bfP_x \left\{k(1-\e^{-f(\xi_t)/k})\right\}
- \bfP_x\bigg\{\int_0^t \beta_k(\xi_s) v^{(k)}_{t-s}(\xi_s) \d s \bigg\} \nonumber \\
& & +\,\bfP_x \bigg\{\int_0^t \bigg[\beta_k(\xi_s)\int_{N(E)}
v^{(k)}_{t-s}(\xi_s)\e^{-\nu(u^{(k)}_{t-s})} F_k(\xi_s,\d\nu)\bigg]\d s\bigg\}   \\
& & +\,\bfP_x \bigg\{\int_0^t \bigg[k\beta_k(\xi_s)\int_{N(E)}
[1-\e^{-\nu(u^{(k)}_{t-s})}] F_k(\xi_s,\d\nu)\bigg]\d s\bigg\},
\nonumber
 \eeqlb
or equivalently
 \beqlb\label{4.4}
v_t^{(k)}(x) +
\int_0^t\bfP_x[\beta_k(\xi_s)\phi_k(\xi_s,v_{t-s}^{(k)}) +
\psi_k(\xi_s,v_{t-s}^{(k)})]\d s = \bfP_x k[1-\e^{-f(\xi_t)/k}],
 \eeqlb
where $\psi_k$ is given by (\ref{3.10}) and (\ref{3.13}), and
 \beqlb\label{4.5}
\phi_k(x,f) = - f(x)\int_{M_0(E)} h_k(x,\pi,1-\pi(f)/k) G(x,\d\pi).
 \eeqlb

\blemma\label{l4.1}
If the conditions of Lemma \ref{l3.1} (i) are fulfilled and if $\beta_k \to \beta\in
B^+(E)$ uniformly as $k\to\infty$, then, for each $a\ge 0$, we have $\phi_k(x,f) \to
f(x)$ uniformly on $E \times B^+_a(E)$ and the solution $v_t^{(k)}(x,f)$ to
(\ref{4.4}) converges boundedly and uniformly on the set $[0,a] \times E\times
B^+_a(E)$ of $(t,x,f)$ to the unique bounded positive solution $V_tf(x)$ to the
evolution equation
 \beqlb\label{4.6}
V_tf(x) - \int_0^t\bigg[\int_E \beta(y)
\zeta(y,V_{t-s}f) P_s(x,\d y)\bigg]\d s = P_t f(x),
 \eeqlb
where $\zeta(\cdot,\cdot)$ is defined by (\ref{3.15}) and (\ref{3.21}).
\elemma

Based on this lemma, one can show as in section 3 that the finite
dimensional distributions of $\{X_t^{(k)}: t\ge0\}$ under
$\bfQ_{(\mu)}^{(k)}$ converge as $k\to\infty$ to those of the
process $\{X_t: t\ge0\}$ with semigroup $(Q_t)_{t\ge0}$ defined by
(\ref{3.22}) and (\ref{4.6}). Since $\alpha_k(x)\equiv 0$ in the
approximating sequence, we call $\{X_t: t\ge0\}$ a {\it rebirth
superprocess}. Note that (\ref{4.6}) is the special form of
(\ref{3.20}) with local branching mechanism $\phi(x,z) \equiv
-\beta(x)z$, which exactly compensates the death factor in the
non-local branching mechanism. This observation might be helpful
in understanding the non-local branching mechanism given by
(\ref{3.21}).

\section{Multitype superprocesses}

\setcounter{equation}{0}

In this section, we deduce the existence of a class of multitype
superprocesses from that of the non-local branching superprocess
constructed in section 3 following the arguments in Li (1992b).
Let $E$ and $I$ be two Lusin topological spaces and let $\xi=
\{\itOmega, (\eta_t, \alpha_t), {\mcr F}, {\mcr F}_t,
\bfP_{(x,a)}\}$ be a right continuous strong Markov process with
state space $E\times I$. Let $\phi (\cdot,\cdot,\cdot)$ and $\zeta
(\cdot,\cdot,\cdot)$ be given by (\ref{3.19}) and (\ref{3.31}),
respectively, with $x\in E$ replaced by $(x,a)\in E\times I$. Let
$\beta(\cdot,\cdot) \in E\times I$ and let $\pi(x,a,\d b)$ be a
probability kernel from $E\times I$ to $I$. As a special form of
the model given in section 3, we have a non-local branching
superprocess $\{X_t(dx,da): t\ge0\}$ in $M(E\times I)$ with
transition probabilities determined by
 \beqlb\label{5.1}
\bfQ_\mu \exp\{-X_t(f)\} = \exp \{-\mu(V_tf)\}, \qquad t\ge0, f\in
B^+(E\times I),
 \eeqlb
where $V_tf$ is the unique bounded positive solution to
 \beqlb\label{5.2}
V_tf(x,a) &+&
\int^t_0\bfP_{(x,a)}[\phi(\eta_s,\alpha_s,V_{t-s}f(\eta_s,\alpha_s))
+ \beta(\eta_s,\alpha_s)V_{t-s}f(\eta_s,\alpha_s)] \d s   \nonumber \\
&-& \int^t_0\bfP_{(x,a)}
[\beta(\eta_s,\alpha_s)\zeta(\eta_s,\alpha_s,
\pi(\eta_s,\alpha_s,V_{t-s}f(\eta_s,\cdot)))] \d s    \\
&=& \bfP_{(x,a)} [f(\eta_t,\alpha_t)].  \nonumber
 \eeqlb
We may call $\{X_t: t\ge0\}$ a {\it multitype superprocess} with type space $I$.
Heuristically,
$\{\eta_t: t\ge0\}$ gives the law of migration of the ``particles'', $\{\alpha_t:
t\ge0\}$ represents the mutation of their types, $\phi(x,a,\cdot)$ describes the
amount of the $a$-type offspring born when an $a$-type parent dies at
$x\in E$, $\zeta(x,a,\cdot)$ describes the amount of the offspring born by this
parent that change into new types randomly according to the kernel $\pi(x,a,\d b)$,
and $\beta(x,a)$ represents the birth rate of the changing-type offspring at
$x\in E$. It is assumed that all of the offspring start migrating from the death
site of their parent. Note that the migration process $\{\eta_t: t\ge0\}$ and
the mutation process $\{\alpha_t: t\ge0\}$ are not necessarily independent.

Now let us consider a special case which has been studied in the literature.
Suppose that $I= \{1,\dots,k\}$ and for each $i\in I$, $\eta^{(i)}$ is a
right continuous strong Markov process in $E$ with semigroup $(P^{(i)}_t)_{t\ge0}$,
$\phi^{(i)}$ belongs to the class given by (\ref{3.19}) and $\zeta^{(i)}$ belongs to
the class given by (\ref{3.31}). Let $\xi$ be a right continuous strong Markov
process in the product space $E\times I$ with transition  semigroup $(P_t)_{t\ge0}$
defined by
 \beqnn
P_tf(x,i)=\int_E f(y,i)P^{(i)}_t(x,\d y),
\qquad f\in B^+(E\times I).
 \eeqnn
Let $\phi((x,i),z) = \phi^{(i)}(x,z)$. Suppose that $\beta^{(i)} \in B^+(E)$ and
$\pi(x,i,\cdot)$ is a Markov kernel from $E\times I$ to $I$ having the decomposition
 \beqnn
\pi(x,i,\cdot) = \sum^k_{j=1}p^{(i)}_j(x)\delta_j(\cdot),
 \eeqnn
where $p_j^{(i)}(x) \ge 0$, $\sum_{j=0}^k p_j^{(i)}(x) \equiv 1$ and
$\delta_j$ denotes the unit mass at $j\in I$. Then we have a multitype
superprocess $\{X_t: t\ge0\}$ in $M(E\times I)$ by (\ref{5.1}) and (\ref{5.2}).
For $i\in I$ and $\mu \in M(E\times I)$ we define $\mu^{(i)}\in M(E)$ by $\mu^{(i)}(B)
= \mu(B\times\{i\})$. The map $\mu \mapsto (\mu^{(1)}, \dots, \mu^{(k)})$
is clearly a homeomorphism between $M(E\times I)$  and the $k$-dimensional
product space $M(E)^k$. Therefore, $\{(X_t^{(1)}, \dots, X_t^{(k)}): t\ge 0\}$ is a
Markov process in the space $M(E)^k$, which may be called a {\it $k$-type
superprocess}. Clearly, this class of $k$-type superprocesses coincides with
the one defined in Li (1992b). Heuristically, $\eta^{(i)}$ gives the law of
the migration of the $i$th type ``particles'', $\phi^{(i)}(x,\cdot)$  describes the
amount of the $i$th type offspring born when an $i$th type parent dies at point $x\in
E$, $\zeta^{(i)}(x,\cdot)$ describes the amount of the offspring born by this parent
that change into new types randomly according to the discrete distribution
$\{p^{(i)}_1(x), \dots, p^{(i)}_k(x)\}$, and $\beta^{(i)}(x)$ represents the birth
rate of the changing-type offspring at $x\in E$. The study of multitype superprocesses
was initiated by Gorostiza and Lopez-Mimbela (1990); see also Gorostiza and Roelly
(1991) and Gorostiza et al (1992).

\section{Superprocess-controlled immigration}

\setcounter{equation}{0}

By the discussions in the last section, we have a special $2$-type superprocess
$\{(X^{(1)}_t, X^{(2)}_t): t\ge 0\}$ in $M(E)^2$ with transition probabilities
determined by
 \beqlb\label{6.1}
\bfQ_{(\mu^{(1)},\mu^{(2)})}
\exp\left\{-X^{(1)}_t(f^{(1)})-X^{(2)}_t(f^{(2)})\right\}
=\exp\left\{-\mu^{(1)}(v^{(1)}_t)-\mu^{(2)}(v^{(2)}_t)\right\},
 \eeqlb
where $v^{(1)}_t(\cdot)$ and $v^{(2)}_t(\cdot)$ are defined uniquely by
 \beqlb\label{6.2}
v^{(1)}_t(x) + \int^t_0\bigg[\int_E
\left(\phi^{(1)}(y,v^{(1)}_{t-s}(y)) - v^{(2)}_{t-s}(y)\right)
P^{(1)}_s(x,\d y)\bigg]\d s = P^{(1)}_t f^{(1)}(x),
 \eeqlb
and
 \beqlb\label{6.3}
v^{(2)}_t(x)
+ \int^t_0\bigg[\int_E \phi^{(2)} (y,v^{(2)}_{t-s}(y)) P^{(1)}_s(x,\d y)\bigg]\d s
= P^{(2)}_t f^{(2)}(x).
 \eeqlb
In particular, if $f^{(2)}\equiv 0$, we have $v^{(2)}_t \equiv 0$ and
 \beqlb\label{6.4}
\bfQ_{(\mu^{(1)},\mu^{(2)})}
\exp\left\{-X^{(1)}_t(f^{(1)})\right\}
=\exp\left\{-\mu^{(1)}(v^{(1)}_t)\right\},
 \eeqlb
where $v^{(1)}_t(\cdot)$ is given by
 \beqlb\label{6.5}
v^{(1)}_t(x) + \int^t_0\bigg[\int_E \phi^{(1)}(y,v^{(1)}_{t-s}(y))
P^{(1)}_s(x,\d y)\bigg]\d s = P^{(1)}_t f^{(1)}(x).
 \eeqlb
Thus $\{X^{(1)}_t: t\ge 0\}$ is a superprocess in $M(E)$ with parameters
$(\eta^{(1)},\phi^{(1)})$. On the other hand, by an expression of weighted
occupation times, the value in (\ref{6.1}) is equal to
 \beqnn
\bfQ_{(\mu^{(1)},\mu^{(2)})}
\exp\left\{-X^{(1)}_t(f^{(1)})\right\}
\exp\left\{-\mu^{(2)}(v^{(2)}_t)-\int_0^tX^{(1)}_s(v^{(2)}_{t-s})\d
s\right\};
 \eeqnn
see e.g. Dawson (1993) and Dynkin (1993, 1994). Then we see that
 {\small\beqnn
\bfQ_{(\mu^{(1)},\mu^{(2)})}
\bigg[\exp\left\{-X^{(2)}_t(f^{(2)})\right\}\bigg|X^{(1)}_s: s\ge
0\bigg] =
\exp\left\{-\mu^{(2)}(v^{(2)}_t)-\int_0^tX^{(1)}_s(v^{(2)}_{t-s})\d
s\right\}.
 \eeqnn$\!\!$}
That is, given $\{X^{(1)}_t: t\ge 0\}$, the second coordinate $\{X^{(2)}_t: t\ge 0\}$
is a superprocess with parameters $(\eta^{(2)},\phi^{(2)})$ and {\it with immigration
controlled by} $\{X^{(1)}_t: t\ge 0\}$. A special class of superprocess-controlled
immigration processes have been studied in Hong and Li (1999) and their relation with
multitype superprocesses has been pointed out in Li (2002).

\section{Mass-structured superprocesses}

\setcounter{equation}{0}

A multitype superprocess $\{X_t(\d x,\d a): t\ge0\}$ with type space $I= (0,\infty)$
can be called a {\it mass-structured superprocess} if we interpret $x\in E$ and
$a>0$ as the coordinates of position and mass, respectively. For the mass-structured
superprocess, we may consider its {\it aggregated process} $\{Y_t: t\ge0\}$
defined by
 \beqlb\label{7.1}
Y_t(\d x) := \int_0^\infty a X_t(\d x,\d a),
\qquad x\in E.
 \eeqlb
Since the integrand on the right hand side is unbounded, $\{Y_t: t\ge0\}$ is only
well-defined under some restrictions. In general, $\{Y_t: t\ge0\}$ is not Markovian.
Let $A$ be the weak generator of $\{(\eta_t,\alpha_t): t\ge0\}$ and
$(T_t)_{t\ge0}$ the locally bounded semigroup of finite kernels on $E\times I$
with generator
 \beqlb\label{7.2}
Bf(x,a)
=
Af(x,a)
+ \beta(x,a)[m(x,a)\pi(a,f(x,\cdot)) - f(x,a)]
- b(x,a)f(x,a),
 \eeqlb
where $b(x,a)$ is the coefficient of the linear term of $\phi(x,a,z)$ and $m(x,a)$
is defined by (\ref{3.32}) with $x\in E$ replaced by $(x,a)\in E\times I$. Indeed,
(\ref{7.2}) is of the same form as (\ref{3.28}) with $(x,a)$ instead of $x$. By the
discussions in section 3, the first moments of the superprocess are given by
 \beqlb\label{7.3}
\bfQ_\mu \{X_t(f)\} = \mu(T_tf), \qquad f\in B^+(E\times I).
 \eeqlb

In practice, we may have that a newborn offspring is no larger than
its parent, which corresponds to the assumption that $\pi(a,\cdot)$ is
supported by $(0,a)$. Let $H(x,a)=a$ and suppose that $H\in {\cal D}(A)$ is a
$c_1$-excessive function of $\{(\eta_t,\alpha_t): t\ge0\}$ for some constant
$c_1>0$. In this case, we have $BH(x,a)\le (c_1+\|b\|)H(x,a)$
and hence $H\in {\cal D}(B)$ is a $(c_1+\|b\|)$-excessive function of $(T_t)_{t\ge 0}$.
It follows from (\ref{7.3}) that
 \beqlb\label{7.4}
\bfQ_\mu \{X_t(H)\} \le \e^{(c_1+\|b\|)t}\mu(H).
 \eeqlb
Then we may change the state space slightly and take any $\sigma$-finite measure
$\mu$ on $E\times (0,\infty)$ satisfying $\mu(H) <\infty$ as the initial state
of $\{X_t: t\ge0\}$; see e.g. El Karoui and Roelly (1991) or Li (1992c).
In this case, (\ref{7.4}) implies that $X_t(H)<\infty$ a.s. for all $t\ge0$ so
that (\ref{7.1}) defines an aggregated process $\{Y_t: t\ge0\}$ with finite
measure values.

A special type of mass-structured superprocess with Markovian
aggregated process has been studied by Gorostiza (1994). Assume
that $\beta(\cdot,\cdot) \equiv 0$ and $\alpha_t = g(t,\alpha_0)$
for a deterministic mapping $g(\cdot,\cdot)$ from
$[0,\infty)\times (0,\infty)$ to $(0,\infty)$. Let $\bfP_x$ denote
the conditional law of $\{\eta_t: t\ge0\}$ given $\eta_0 =x$. For
$f\in B^+(E)$, (\ref{5.2}) becomes
 \beqlb\label{7.5}
V_tf(x,a) + \int^t_0\bfP_x [\phi(\eta_s,g(s,a),
V_{t-s}f(\eta_s,g(s,a)))] \d s = \bfP_x f(\eta_t).
 \eeqlb
Since the motion of $\alpha_t = g(t,\alpha_0)$ is deterministic,
if $X_0$ is supported by $E \times \{a\}$, then $X_t$ is supported
by $E \times \{g(t,a)\}$ and $Y_t = g(t,a) X_t$. In this case,
$\{Y_t: t\ge0\}$ is a Markov process since the transformation $X_t
\mapsto Y_t$ loses no information. For $B\in {\mcr B}(E)$, let
$X^a_t(B) = X_t(B \times \{g(t,a)\})$. Then $\{X^a_t: t\ge 0\}$ is
an inhomogeneous superprocess with cumulant semigroup $(V^a_{r,t})
_{t\ge r\ge0}$ defined by $V^a_{r,t}f(x) := V_{t-r}f(x,g(r,a))$,
which has underlying process $\{\eta_t: t\ge0\}$ and
time-dependent branching mechanism $\phi (x,g(t,a),\cdot)$. This
gives a representation of the aggregated process in terms of an
inhomogeneous superprocess. A representation of this type was
first given by Gorostiza (1994) in the case where $\alpha_t =
\alpha_0\e^{ct}$ for a constant $c\in \dfR$. Gorostiza (1994)
obtained the process as high density limit of a sequence of
branching particle systems where the mass of each offspring is
equal to that of its parent multiplied by a fixed positive
constant factor, and the mass of any particle does not change
during its lifetime, realizing in a particular case the program
mentioned at the end of section 3.

\section{Multilevel superprocesses}

\setcounter{equation}{0}

Multilevel superprocesses arise as limits of multilevel branching
particle systems. In a two level system, objects at the higher
level consist of non-trivial sub-populations of objects at  the
lower level and both lower level and higher level objects can
branch. A lower level object consisting of a population
can be described by a measure on some space $S$. We can
then view a two level system as a multitype system with
$I=M(S)^\circ$, the space of non-trivial finite Borel measures on
$S$. Non-local branching is natural in this context. For example,
at the particle level the offspring of a second order object
consisting of a set of particles could consist of a subset of the
particles or include more than one copy of the original particles.

To make this precise, we may let $S$ be a topological Lusin space and
$\{\alpha_t: t\ge0\}$ be the Markov process with state space
$M(S)^\circ$ obtained by killing a superprocess at its extinction
time. Then $\{X_t:t\geq 0\}$ is a Markov process with
state space $M(E\times M(S)^\circ)$, which can be called a {\it
multilevel superprocess} generalizing the model of Dawson and
Hochberg (1991), Dawson et al (1990) and Wu (1994).
For the multilevel process, it is also natural to study the {\it
aggregated process} $\{Y_t: t\ge0\}$ defined by
 \beqlb\label{8.1}
Y_t(A\times B) := \int_A\int_{M(S)} a(B) X_t(\d x,\d a),
\qquad A\in {\cal B}(E), B\in {\cal B}(S).
 \eeqlb

To illustrate the possibilities of non-local branching, consider
the case in which $E$ is a singleton. In this case, we may view $\{X_t: t\ge0\}$
as a superprocess with state space $M(M(S)^\circ)$. A possible non-local branching
mechanisms is obtained by taking $\pi (\mu,\d\nu )$ to be the law of
 \beqlb\label{8.2}
\frac{\mu (S)}{N}\sum_{i=1}^{N}\delta _{Z_i},
 \eeqlb
where $N\ge1$ is an integer-valued random variable, and $\{Z_1,Z_2,\cdots\}$ are
i.i.d. random variables in $S$ with distribution $\mu(S)^{-1}\mu(\cdot)$. That
is, the offspring of a level two object $\mu\in M(S)$ is a single point measure
with the same total mass as $\mu$ and its location is selected randomly according
to the empirical measure of a sample from the normalized parent distribution.

Another possibility is
given by
 \beqlb\label{8.3}
\pi (\mu,\d\nu) = \delta_{\mu_B}(\d\nu),
 \eeqlb
where $B\in {\cal B}( S)$ and $\mu_B\in M(S)$ is defined by $\mu_B(A) =
\mu(A\cap B)$. In this case the offspring is a level two object in which only
level one individuals falling in the set $B\subset S$ are present.

In the case in which $E$ is a countable set, we may interpret the $M(E\times
M(S)^\circ)$-valued process $\{X_t: t\geq 0\}$ as a population in a sequence of
islands. The $E$-coordinate tells in which island the $M(S)^\circ$-valued objects
$\{\alpha _t: t\ge0\}$ in the first level are located. The non-local branching
is given by $\pi (x,\mu ,d\nu ) = \pi(\mu ,d\nu)$, which only acts on the
$M(S)^\circ$-coordinate at the higher level. Jumps in $E$ described by
$\{\eta_t: t\ge0\}$ correspond to the independent migration of ``clans''
(families of the lower level) between the islands. Suggestively, we may
call $\{X_t: t\geq 0\}$ a {\it stepping stone type superprocess}.

Properties and applications of multilevel superprocesses involving local
branching have been studied extensively in the literature; see e.g. Dawson and
Hochberg (1991), Dawson et al (1990, 1994, 1995), Etheridge (1993), Gorostiza
(1996), Gorostiza et al (1995), Hochberg (1995), Wu(1994) and the references
therein.

\section{Age-reproduction-structured superprocesses}

\setcounter{equation}{0}

Let $E$ be a Lusin topological space and let $\xi= \{\itOmega,
(\eta_t, \alpha_t, \theta_t), {\mcr F}, {\mcr F}_t,
\bfP_{(y,a,z)}, \gamma\}$ be a Borel Markov process with state
space $E\times \dfR^+ \times \dfN^+$, where $\gamma$ is a terminal
time; see Sharpe (1988, p.65). We assume that both $\alpha_t$ and
$\theta_t$ are non-decreasing processes. Let
$\beta(\cdot,\cdot,\cdot) \in B(E\times \dfR^+ \times \dfN^+)$ and
let $\zeta (\cdot,\cdot,\cdot,\cdot)$ be given by (\ref{3.31})
with $x\in E$ replaced by $(x,a,z) \in E\times \dfR^+ \times
\dfN^+$. As a special form of the models given in sections 4 and
5, we have a rebirth multitype superprocess $\{X_t: t\ge0\}$ in
$M(E\times \dfR^+ \times \dfN^+)$ with transition probabilities
determined by
 \beqlb\label{9.1}
\bfQ_\mu \exp\{-X_t(f)\} = \exp \{-\mu(V_tf)\}, \qquad t\ge0, f\in
B^+(E\times \dfR^+ \times \dfN^+),
 \eeqlb
where $V_tf$ is the unique bounded positive solution to
 \beqlb\label{9.2}
V_tf(y,a,z) &-& \int^t_0\bfP_{(y,a,z)}
[\beta(\eta_s,\alpha_s,\theta_s)\zeta(\eta_s,\alpha_s,\theta_s,
V_{t-s}f(\eta_s,0,0)) 1_{\{\alpha_s < \gamma\}}] \d s    \\
&=& \bfP_{(y,a,z)} [f(\eta_t,\alpha_t,\theta_t)1_{\{\alpha_s <
\gamma\}}].  \nonumber
 \eeqlb
It is not hard to check that the first moments of the superprocess  are given by
 \beqlb\label{9.3}
\bfQ_\mu \{X_t(f)\} = \mu(T_tf), \qquad f\in B^+(E\times \dfR^+
\times \dfN^+),
 \eeqlb
where $(T_t)_{t\ge 0}$ is a semigroup of bounded linear operators
on $B^+(E\times \dfR^+ \times \dfN^+)$ defined by
 \beqlb\label{9.4}
T_tf(y,a,z) &-& \int^t_0\bfP_{(y,a,z)}
[\beta(\eta_s,\alpha_s,\theta_s)m(\eta_s,\alpha_s,\theta_s)
T_{t-s}f(\eta_s,0,0) 1_{\{\alpha_s < \gamma\}}] \d s    \\
&=& \bfP_{(y,a,z)} [f(\eta_t,\alpha_t,\theta_t) 1_{\{\alpha_s <
\gamma\}}],  \nonumber
 \eeqlb
where $m(\cdot,\cdot,\cdot)$ is given by (\ref{3.32}) with $x\in
E$ replaced by $(x,a,z) \in E\times \dfR^+ \times \dfN^+$. Using
$(T_t)_{t\ge 0}$ we may rewrite (\ref{9.2}) into
 \beqlb\label{9.5}
V_tf(y,a,z) + \int^t_0T_s[\beta mV_{t-s}f
- \zeta(\cdot,\cdot,\cdot, V_{t-s}f(\cdot,0,0))] (y,a,z)\d s
= T_tf(y,a,z).
 \eeqlb

In the case $\alpha_t \equiv \alpha_0 + t$, we may call $\{X_t: t\ge0\}$ an
{\it age-reproduction-structured superprocess}. Heuristically, $\eta_t$ represents
the location of a ``particle'', $\alpha_t$ its age and $\theta_t$ the number of its
offspring born in the time interval $(t-\alpha_t,t]$. At each branching time,
the particle gives birth to a random number of offspring whose motions start
from the branching site and whose ages and reproduction numbers start from zero.
The particle does not disappear at its branching times, it is removed from the
population only when its age exceeds the lifetime $\gamma$. An interesting limit
theorem for age-reproduction-structured branching particle systems was proved
in Bose and Kaj (2000) which leads to the superprocess in the special case where
$E$ is a singleton and $\eta_t \equiv \eta_0$. (Compare (\ref{9.5}) and their
equation (2.8).)

\end{document}